\documentclass [10pt,twoside,b5paper]{article}
\usepackage{amsfonts}
\usepackage{amsthm}
\usepackage{amsmath}
\usepackage{amstext}
\usepackage{amssymb}
\usepackage{mathrsfs}
\usepackage{epsf}              
\usepackage{graphicx}          
\usepackage{fancybox}          
\usepackage{color}             
\usepackage{fancyhdr}
\usepackage[hang,footnotesize]{caption2}  

\def\mathcal{\mathscr}
\newfont{\aaa}{cmb10 at 18pt}
\newfont{\bbb}{cmb10 at 10pt}
\newtheorem{lem}{Lemma}[]
\newtheorem{defi}{Definition}[]
\newtheorem{thm}{Theorem}
\newtheorem{cor}{Corollary}[]
\newtheorem{rem}{Remark}[]
\newtheorem{pro}{Proposition}[]
\def\QED{\hfill$\Box$}

\def\LL{{\mathcal L}}

\def\Z{\mathbb Z}

\def\dis{\displaystyle}
\newcommand{\beq}{\begin{equation}}
\newcommand{\eeq}{\end{equation}}
\newcommand{\bey}{\begin{eqnarray}}
\newcommand{\eey}{\end{eqnarray}}
\newcommand{\beyy}{\begin{eqnarray*}}
\newcommand{\eeyy}{\end{eqnarray*}}

\renewcommand{\headsep}{0.7cm}
 \makeatletter
\def\@evenhead{
   \vbox{\hbox to \textwidth
{}{\hspace{0mm}{\footnotesize \thepage}}{\hspace{7.5cm}
   {\footnotesize {Yucai SU, Lamei YUAN}}}
  \protect\vspace{1truemm}\relax
   \hrule depth0pt height0.15truemm width\textwidth
   }}
   \def\@evenfoot{}
\def\@oddhead{
    \vbox{\hbox to \textwidth
  {{\hspace{0cm}{\footnotesize Quantization of the Schr\"{o}dinger-Virasoro Lie algebra}
  \hfill{\footnotesize \thepage}}\hspace{0mm}}{}
   \protect\vspace{1truemm}\relax
   \hrule depth0pt height0.15truemm width\textwidth
  }}
  \def\@oddfoot{}
\makeatother

\begin{document}

\thispagestyle{empty} \thispagestyle{fancy} {
  \fancyhead[lO,RE]{\footnotesize  Front.  Math. China \\
  DOI 10.1007/s11464-010-0036-2\\[3mm]
  \includegraphics[0,-50][0,0]{11.bmp}}
  \fancyhead[RO,LE]{\scriptsize \bf 
} \fancyfoot[CE,CO]{}}
\renewcommand{\headrulewidth}{0pt}
\renewcommand{\headsep}{0.7cm}


\setcounter{page}{1}
\qquad\\[8mm]

\noindent{\aaa{Quantization of the Schr\"{o}dinger-Virasoro Lie
algebra
$^{^{^{\displaystyle*}}}$}}\\[1mm]

\noindent{\bbb Yucai SU,\quad Lamei YUAN}\\[-1mm]

\noindent\footnotesize{Department of Mathematics, University of
Science\! and \!Technology \!of\! China, Hefei 230026,
China}\\[6mm]

 \vskip-2mm
 \noindent{\footnotesize$\copyright$ Higher Education Press and
Springer-Verlag Berlin Heidelberg 2010}
 \vskip 4mm

\normalsize\noindent{\bbb Abstract}\quad In this paper, we use the
general quantization method by Drinfel'd twists to quantize the
Schr\"{o}dinger-Virasoro Lie algebra whose Lie bialgebra structures
were recently discovered by Han-Li-Su. We give two different kinds
of Drinfel'd twists, which are then used to construct the
corresponding
 Hopf algebraic structrues. Our results extend the class of
examples of noncommutative and noncocommutative Hopf algebras.
\vspace{0.3cm}

\noindent{\bbb Keywords}\quad Lie bialgebras, quantization,
Schr\"{o}dinger-Virasoro Lie algebra\\
{\bbb MSC}\quad 17B05, 17B37,
17B62, 17B68\\[0.4cm]

\noindent{\bbb{1\quad Introduction}}\\[0.1cm]
In Hopf algebra or quantum group theory, two standard methods to
yield new bialgebras from old ones are by twisting the product by a
2-cocycle but keeping the coproduct unchanged, and by twisting the
coproduct by a Drinfel'd twist but preserving the product.
Constructing quantizations of Lie bialgebras is an important
approach to producing new quantum groups (see \cite{D1, D2, ES} and
references therein). Drinfel'd in \cite{D3} formulated a number of
problems in quantum group theory, including the existence of a
quantization for Lie bialgebras. In the paper \cite{EK1} Etingof and
Kazhdan gave a positive answer to some of Drinfeld's questions. In
particular, they showed the existence of quantizations for Lie
bialgebras, namely, any classical Yang-Baxter algebra can be
quantized. Since then the interests in quantizations of Lie
bialgebras have been growing in the mathematical literatures (e.g.,
\cite{G, EH, EK6, SS, HW}). \par The Schr\"{o}dinger-Virasoro Lie
algebra considered in this paper was introduced in the context of
non-equilibrium statistical physics during the process of
investigating the free
Schr\"{o}dinger equations (see \cite{H1, H3}). This Lie algebra  \\

\vspace{0mm}
\noindent \hrulefill\hspace{117mm}\\
{\footnotesize $^*$  Received November 29, 2009; accepted June 10,
2010}

\vspace{-2mm}
\noindent \hrulefill\hspace{117mm}\\
{\footnotesize Corresponding author: Lamei YUAN, E-mail:
lmyuan@mail.ustc.edu.cn} \\
is closely related to Schr\"{o}dinger algebra and Virasoro algebra,
both of which play important roles in many areas of mathematics and
physics (e.g., statistical physics). The Schr\"{o}dinger-Virasoro
Lie algebra, denoted by $\LL$, is an infinite-dimensional vector
space with basis $\{L_n,Y_p,M_n\,\big|\,n\in \Z,p\in
\frac{1}{2}+\Z\}$ and the following non-vanishing Lie brackets
\begin{eqnarray}\label{LB}\begin{array}{lll}
&&[L_m,L_{n}]=(n-m)L_{n+m},\ \ \
[L_m,M_n]=nM_{n+m},\\[6pt]
&&[\,L_n,Y_p\,]=(p-\frac{n}{2})Y_{p+n},\ \ \ \ \,\
[\,\,Y_p\,,Y_{q}\,\,]=(q-p)M_{p+q},\end{array}
\end{eqnarray}
for all $m,n\in\Z$ and $p,q\in\Z+1/2.$ This kind of Lie algebras has
been investigated in a number of
 papers. Some of these investigations \cite{GJP, LS2, LS3, LSZ} focus on its structure
 theory including derivations, central extension and automorphism groups, others \cite{LS1, RU, TZ, U} on its representations.
\par Recently, the
Lie bialgebra structures on the Schr\"{o}dinger-Virasoro Lie algebra
$\LL$ were discussed in \cite{HS}, which turned out to be not all
coboundary triangular (for definition, see p.28, \cite{ES}). In the
present paper, we use the general quantization method by Drinfel'd
twists (cf.~\cite{D1}) to quantize explicitly the newly determined
triangular Lie bialgebra structures on the Schr\"{o}dinger-Virasoro
Lie algebra $\LL$. Actually, this process completely depends on the
construction of Drinfel'd twists determined by the $r$-matrix
(namely, the triangular Lie bialgebra structures on the
Schr\"{o}dinger-Virasoro Lie algebra). Our results extend the class
of examples of noncommutative and noncocommutative Hopf
algebras.\par
 The main results of this paper are listed as follows:
\begin{thm}\label{theo1}
With the choice of two distinguished elements
$h:=\frac{1}{n_0}L_{0}$ and $e:=M_{n_0}$ $(n_0\ne0)$ such that
$[h,e]=e$ in  $\mathcal {L}$, there exists a structure of
noncommutative and noncocommutative Hopf algebra $(U(\mathcal
{L})[[t]], m,\iota,\Delta, S,\epsilon)$ on $U(\mathcal {L})[[t]]$
over $\mathbb{F}[[t]]$ with $U(\mathcal {L})[[t]]/tU(\mathcal
{L})[[t]]\cong U(\LL)$, which leaves the product and counit of
$U(\mathcal {L})[[t]]$ undeformed but with the deformed
comultiplication and antipode defined by:
\begin{eqnarray*}
\Delta{(L_n)}&:=&1\otimes
L_n+L_n\otimes(1-et)^{\frac{n}{n_0}}+n_0h\otimes(1-et)^{-1}M_{n+n_0}t,\\[1mm]
\Delta{(M_k)}&:=&1\otimes M_k+M_k\otimes (1-et)^{\frac{k}{n_0}},\\[1mm]
\Delta{(Y_p)}&:=&1\otimes Y_p+Y_p\otimes(1-et)^{\frac{p}{n_0}},\\[1mm]
S(L_n)&:=&-(1-et)^{-\frac{n}{n_0}}(L_n-n_0M_{n+n_0}h^{[1]}_1t),\\[1mm]
S(M_k)&:=&-(1-et)^{-\frac{k}{n_0}}\cdot M_k,\\[1mm]
S(Y_p)&:=&-(1-et)^{-\frac{p}{n_0}}\cdot Y_p.
\end{eqnarray*}
\end{thm}
\begin{thm}\label{theo2}\vskip-3pt
With the choice of two distinguished elements
$h:=\frac{2}{n_0}L_{0}$ and $e:=Y_{\frac{n_0}{2}}$ $(n_0\in2\Z+1)$
such that $[h,e]=e$ in $\mathcal {L}$, there exists a structure of
noncommutative and noncocommutative Hopf algebra $(U(\mathcal
{L})[[t]], m,\iota, \Delta, S,\epsilon)$ on $U(\mathcal {L})[[t]]$
over $\mathbb{F}[[t]]$ with $U(\mathcal {L})[[t]]/tU(\mathcal
{L})[[t]]\cong U(\LL)$, which leaves the product and counit of
$U(\mathcal {L})[[t]]$ undeformed but with the deformed
comultiplication and antipode defined by:
\begin{eqnarray*}
\Delta({L_n})&:=&1\otimes L_n+L_n\otimes (1-et)^{\frac{2n}{n_0}}
+\frac{n_0-n}{2}h\otimes (1-et)^{-1}Y_{n+\frac{n_0}{2}}t\\[1mm]
&&+\frac{n(n-n_0)}{4}h^{(2)}\otimes (1-et)^{-2}M_{n+n_0}t^2,\\[1mm]
\Delta({M_k})&:=&1\otimes M_k+M_k\otimes (1-et)^{\frac{2k}{n_0}},\\[1mm]
\Delta({Y_p})&:=&1\otimes
Y_p+Y_p\otimes(1-et)^{\frac{2p}{n_0}}-(p-\frac{n_0}{2})h\otimes
(1-et)^{-1}M_{p+\frac{n_0}{2}}t,\\[1mm]
S(L_n)&:=&-(1-et)^{-\frac{2n}{n_0}}\big(L_n+\frac{n-n_0}{2}Y_{n+\frac{n_0}{2}}h^{[1]}_1
t+\frac{n(n-n_0)}{4}M_{n+n_0}h^{[2]}_{2}t^2\big),\\[1mm]
S(Y_p)&:=&-(1-et)^{-\frac{2p}{n_0}}\big(Y_p+(p-\frac{n_0}{2})M_{p+\frac{n_0}{2}}h^{[1]}_1t\big),\\[1mm]
S(M_k)&:=&-(1-et)^{-\frac{2k}{n_0}}M_k.
\end{eqnarray*}
\end{thm}\par
 Throughout this paper $\mathbb{F}$
denotes a field of characteristic zero. All vector spaces and
algebras are assumed to be over $\mathbb{F}$. $\Z$, $\Z_+$ and
$\Z^*$ stand for the sets of integers, nonnegative and nonzero
integers respectively.
\\[4mm]

\noindent{\bbb 2\quad Preliminaries}\\[1mm]

\noindent In this section, we summarize some basic definitions and
results concerning Lie bialgebra structures which will be used in
the following discussions. For a detailed discussion of this subject
we refer the reader to the literatures (e.g. \cite{SS, HS} and
references therein). \par Let $\LL$ be the Schr\"{o}dinger-Virasoro
Lie algebra defined in (\ref{LB}) and $U(\LL)$ the universal
enveloping algebra of $\LL$. Then $U(\mathcal {L})$ is equiped with
 a natural Hopf algebraic structure $(U(\mathcal {L}),m,\iota,\Delta_{0},S_{0},\epsilon)$, i.e.,
\begin{eqnarray}\label{hop}
\Delta_{0}(X)=X\otimes 1+1\otimes X,& S_{0}(X)=-X,
&\epsilon(X)=0\mbox{ \ for }X\in{\mathcal{L}}.
\end{eqnarray}
where $\Delta_{0}$ is a comultiplication, $\epsilon$ is a counit and
$S_0$ is an antipode. In particular, $$\Delta_0(1)=1\otimes 1\ \
\mbox{and} \ \ \epsilon(1)=S_0(1)=1.$$\par In order to search for
the solutions of the Yang-Baxter quantum equation, Drinfel'd in
\cite{D1} introduced the notion of Lie bialgebras in 1983. Since
then, a great deal of attention has been paid to the study of the
quantization of Lie bialgebras as well as Lie bialgebra structures
of some Lie algebras (e.g., \cite{EK1, EK6, GZ, HS}). The following
result is due to \cite{SS}.
\begin{pro}\label{p}
Let $L$ be a Lie algebra containing two linear independent elements
$a,b$ satisfying $[a,b]=kb$ with $0\neq k\in\mathbb{F}.$ Set
$r=a\otimes b-b\otimes a$ and define a linear map
\begin{eqnarray}\label{cybe}\Delta_{r}(x)=x\cdot r=[x,a]\otimes b-b\otimes
[x,a]+a\otimes[x,b]-[x,b]\otimes a,\ \ \mbox{for}\ \ x\in L.
\end{eqnarray} Then $\Delta_r$ equips $L$ with a structure of triangular
coboundary Lie bialgebra.
\end{pro}\par
Equation (\ref{cybe}) implies that $\Delta_r$ is an inner derivation
of $L$. For the Schr\"{o}dinger-Virasoro Lie algebra $\LL$ defined
in (\ref{LB}), it is shown in \cite{HS} that a Lie bialgebra
$(\LL,[\cdot,\cdot],\Delta)$ is triangular coboundary if and only if
$\Delta$ is an inner derivation, which is determined by the
classical Yang-Baxter $r$-matrix $r$. From the above proposition, we
notice that the classical Yang-Baxter $r$-matrix $r$ is uniquely
expressed as the antisymmetric tensor of two distinguished elements
$a,b$ up to nonzero scalars satisfying $[a,b]=kb$ ($k\neq 0$). In
fact, for a given $r$-matrix, we may take two distinguished elements
of the form $h:=k^{-1}a$ and $e:=kb$ such that $[h,e]=e$ with $0\neq
k\in \mathbb{F}$.
\begin{defi}\label{def2}\rm
Let $(H,m,\iota,\Delta_{0},S_{0},\epsilon)$ be a Hopf algebra over a
commutative ring $R$. A Drinfel'd twist $\mathcal {F}$ on $H$ is an
invertible element of $H\otimes H$ such that
$$(\mathcal {F}\otimes 1)(\Delta_{0}\otimes {\rm Id})(\mathcal {F})=(1\otimes \mathcal {F})({\rm Id}\otimes\Delta_{0})(\mathcal
{F}),$$
$$(\epsilon\otimes {\rm Id})(\mathcal {F})=1\otimes 1=({\rm Id}\otimes\epsilon)(\mathcal
{F}).$$
\end{defi}\par
The following result is well known (see \cite{D1, ES}, etc.).
\begin{lem}\label{Legr}\rm
Let $(H,m,\iota,\Delta_{0},S_{0},\epsilon)$ be a Hopf algebra over a
commutative ring and  $\mathcal {F}$ a Drinfel'd twist on $H$, then
$w=m({\rm Id}\otimes S_{0})(\mathcal {F})$ is invertible in $H$ with
$w^{-1}=m(S_{0}\otimes {\rm Id})(\mathcal {F}^{-1})$.  Moreover, if
we define $\Delta$ :$H\rightarrow H\otimes H$ and $S$ :$H\rightarrow
H$ by
$$\begin{array}{llll}
\Delta(x)=\mathcal {F}\Delta_{0}(x)\mathcal {F}^{-1},&
S=wS_{0}(x)w^{-1},  \mbox{ \ for all \ }x\in H.
\!\!\!\!\!\!\!\!\!\!\!\!\!\!\!\!\!\!\!\!\!\!\!\!\end{array}$$ Then
 $(H,m,\iota,\Delta,S,\epsilon)$ is a new Hopf algebra, which is called the
 {\it twisting} of $H$ by the Drinfel'd twist $\mathcal {F}$.
\end{lem}\par Let $\mathbb{F}[[t]]$ be a ring of formal power
series. Assume that $L$ is a triangular Lie bialgebra with a
classical Yang-Baxter $r$-matrix $r$ (see \cite{D1,ES}). Denote by
$U(L)$ the universal enveloping algebra of $L$, with the standard
Hopf algebra structure $(U( L),m,\iota,\Delta_{0},S_{0},\epsilon)$.
Now let us consider the topologically free $F[[t]]$-algebra
$U(L)[[t]]$ (for definition, see p.$4$, \cite{ES}), which can be
viewed as an associative $\mathbb{F}$-algebra of formal power series
with coefficients in $U(L)$. Naturally, $U(L)[[t]]$ is equiped with
an induced Hopf algebra structure arising from that on $U(L)$. By
abuse of notation, we denote it by $(U(
L)[[t]],m,\iota,\Delta_{0},S_{0},\epsilon)$.
\begin{defi}\rm (See Definition 1.4, \cite{HW}) For a triangular Lie bialgebra $L$, $U(L)[[t]]$ is
called a quantization of $U(L)$ by a Drinfel'd twist $\mathcal {F}$
over $U(L)[[t]]$ if $U(L)[[t]]/tU(L)[[t]]\cong U(L)$ and $\mathcal
{F}$ is determined by its $r$-matrix $r$ (namely, its Lie bialgebra
structure).
\end{defi}\par

An algebra $A$ equipped with a classical Yang-Baxter $r$-matrix $r$
is called a {\it classical Yang-Baxter algebra}. It is showed in
\cite{EK1} that any classical Yang-Baxter algebra can be
quantized.\par
 For any element $x$
of a unital $R$-algebra ($R$ a ring) and $a\in R$, we set (see,
e.g., \cite {GZ})
\begin{eqnarray*}
&&x^{(n)}_{a}:=(x+a)(x+a+1)\cdots(x+a+n-1)\\
&&x^{[n]}_{a}:=(x+a)(x+a-1)\cdots(x+a-n+1)
\end{eqnarray*}
and $x^{(n)}:=x^{(n)}_{0}$,  $x^{[n]}:=x^{[n]}_{0}$.

\begin{lem}\label{Legr3}(See \cite{GZ, G}) For any element $x$ of a
unital $\mathbb{F}$-algebra, $a,b\in \mathbb{F}$, and  $r,s,t\in
\mathbb{Z}$, one has
\begin{eqnarray}\label{fir-e}
&&x^{(s+t)}_{a}=x^{(s)}_{a}x^{(t)}_{a+s},\ \ \ \ \  x^{[s+t]}_{a}=
x^{[s]}_{a}x^{[t]}_{a-s},\ \ \ \ \
 x^{[s]}_{a}= x^{(s)}_{a-s+1},\label{fir-e}\\
&&\mbox{$\sum\limits_{s+t=r}$}\frac{(-1)^{t}}{s!t!}x^{[s]}_{a}x^{(t)}_{b}={a-b\choose
r}=
\frac{(a-b)\cdots(a-b-r+1)}{r!},\label{fir-e1}\\
&&\mbox{$\sum\limits_
{s+t=r}$}{\displaystyle\frac{(-1)^{t}}{s!t!}}x^{[s]}_{a}x^{[t]}_{b-s}={\displaystyle{a-b+r-1\choose
r}}=\frac{(a-b)\cdots(a-b+r-1)}{r!}.\label{fir-e2}
\end{eqnarray}
\end{lem}
\begin{rem}\rm One can see that the right-hand sides
of the last two equations do not depend on $x$ from the proof
process (see e.g., Lemma 3, \cite{G}).
\end{rem}
\par
The following popular result will be frequently used in the third
part of this paper.

\begin{lem}\rm \label{Legr4}(see e.g., Proposition 1.3(4), \cite{SF}) For any elements $x, y$ of an
associative algebra $A$, and $m\in\mathbb{Z}_{+}$, one has
\begin{eqnarray}
xy^{m}=\mbox{$\sum\limits_{k=0}^{m}$}(-1)^k{m\choose k}y^{m-k}({\rm
ad\,}y)^k(x).
\end{eqnarray}
\end{lem}
\vspace{4mm}

\noindent{\bbb 3\quad Proof of the main
results}\\[1mm]

\noindent To describe quantizations of $U(\mathcal {L})$, we need to
construct explicitly Drinfel'd twists according to Lemma \ref{Legr}.
Fix $n_0\in \mathbb{Z}^*$. Set
\begin{eqnarray*}
h:=\frac{1}{n_0}L_0,\ \ e:=M_{n_0},\mbox{ \ \ or \ }
h:=\frac{2}{n_0}L_0,\ \ e:=Y_{\frac{n_0}{2}}\mbox{ \ (only when
$n_0$ is odd)}.
\end{eqnarray*}
Clearly, one has $[h,e]=e$ by equation (\ref{LB}). Then from
Proposition \ref{p} it follows that the Lie bialgebra $(\LL,
[\cdot,\cdot], \Delta_r)$ with $r=h\otimes e-e\otimes h$ is
triangular. We shall quantize this triangular Lie bialgebra
structure in this section. To do this, we need some necessary
calculations, which are useful to the construction of Drinfel'd
twists in the sequel.
\begin{lem}\rm \label{lemm1}
For $a\in\mathbb{F}$, $i\in\mathbb{Z}_{+}$, $n\in\mathbb{Z}$, and
$p\in\mathbb{Z}+\frac{1}{2}$, the following  hold in
$U(\mathcal{L})$:\begin{itemize}\parskip-2pt

\item[\rm(i)] If $h=\frac{1}{n_0}L_0$, and $e=M_{n_0}$, then
\begin{eqnarray}
&L_nh^{(i)}_a=h^{(i)}_{a-\frac{n}{n_0}}L_n,&
L_nh^{[i]}_a=h^{[i]}_{a-\frac{n}{n_0}}L_n,\label{fir-eeee}
\\
&M_nh^{(i)}_a=h^{(i)}_{a-\frac{n}{n_0}}M_n,&M_nh^{[i]}_a=h^{[i]}_{a-\frac{n}{n_0}}M_n,
\nonumber\\
&Y_ph^{(i)}_a=h^{(i)}_{a-\frac{p}{n_0}}Y_P,&
Y_ph^{[i]}_a=h^{[i]}_{a-\frac{p}{n_0}}Y_p.\nonumber
\end{eqnarray}
\item[\rm(ii)] If $h=\frac{2}{n_0}L_0$, and $e=Y_{\frac{n_0}{2}}$,
then
\begin{eqnarray*}
&L_nh^{(i)}_a=h^{(i)}_{a-\frac{2n}{n_0}}L_n,&
L_nh^{[i]}_a=h^{[i]}_{a-\frac{2n}{n_0}}L_n,
\\
&M_nh^{(i)}_a=h^{(i)}_{a-\frac{2n}{n_0}}M_n,&M_nh^{[i]}_a=h^{[i]}_{a-\frac{2n}{n_0}}M_n,
\\&Y_ph^{(i)}_a=h^{(i)}_{a-\frac{2p}{n_0}}Y_P,&
Y_ph^{[i]}_a=h^{[i]}_{a-\frac{2p}{n_0}}Y_p.
\end{eqnarray*}

\item[\rm(iii)] In both cases,
\begin{eqnarray}\label{fir-eeee2}
e^nh^{(i)}_{a}=h^{(i)}_{a-n}e^n,&& e^nh^{[i]}_{a}=h^{[i]}_{a-n}e^n.
\end{eqnarray}\end{itemize}\end{lem} \par
\noindent{\it Proof.~}~We only prove the first equation of
(\ref{fir-eeee}) (the others   can be obtained similarly). We have
$[L_n, h]=-\frac{n}{n_0}L_{n}$ and $L_nh=(h-\frac{n}{n_0})L_n$,
i.e., it holds  for $i=1$. Suppose that it holds for $i$, then we
have
\begin{eqnarray*}
L_nh^{(i+1)}_a&=&L_nh^{(i)}_a(h+a+i)=h^{(i)}_{a-\frac{n}{n_0}}L_n(h+a+i)\\
&=&h^{(i)}_{a-\frac{n}{n_0}}(h-\frac{n}{n_0}+a+i)L_n=h^{(i+1)}_{a-\frac{n}{n_0}}L_n.
\end{eqnarray*}
\QED

For any $a\in\mathbb{F}$, we set
\begin{eqnarray}\label{fom1}
&\mathcal
{F}_a=\sum\limits^{\infty}_{r=0}\frac{(-1)^r}{r!}h^{[r]}_a\otimes
e^r t^r,& F_a=
\mbox{$\sum\limits^{\infty}_{r=0}$}\frac{1}{r!}h^{(r)}_a\otimes e^r
t^r,
\\
&u_a=m\cdot (S_0\otimes{ \rm Id})(F_a), &v_a=m\cdot({\rm Id}\otimes
S_0)(\mathcal {F}_a).\nonumber
\end{eqnarray}
Write $\mathcal {F}=\mathcal {F}_0$, $F=F_0$, $u=u_0$, $v=v_0$.
Since $S_0(h^{(r)}_a)=(-1)^rh^{[r]}_{-a}$ and $S_0(e^r)=(-1)^r e^r$,
we have
\begin{eqnarray}\label{fom2}
u_a=\mbox{$\sum\limits^{\infty}_{r=0}$}\frac{(-1)^r}{r!}h^{[r]}_{-a}e^r
t^r,&&v_a= \mbox{$\sum\limits^{\infty}_{r=0}$}\frac{1}{r!}h^{[r]}_a
e^r t^r.
\end{eqnarray}

\begin{lem}\rm \label{lemm2}
For any $a,b\in\mathbb{F}$, one has
\begin{eqnarray*}
\mathcal {F}_aF_b=1\otimes(1-et)^{a-b},&& v_au_b=(1-et)^{-(a+b)}.
\end{eqnarray*}
\end{lem}
\noindent{\it Proof.}~~By  equations (\ref{fir-e1}) and
(\ref{fom1}), we have
\begin{eqnarray*}
\mathcal {F}_aF_b&=&
\mbox{$\sum\limits^{\infty}_{r,s=0}$}\frac{(-1)^r}{r!s!}h^{[r]}_ah^{(s)}_b\otimes
e^r e^s t^r t^s\\&=&\mbox{$\sum\limits^{\infty}_{m=0}$}(-1)^m
\big(\mbox{$\sum\limits_{r+s=m}$}\frac{(-1)^r}{r!s!}h^{[r]}_{a}h^{(s)}_b\big)\otimes e^m t^m\\
&=&\mbox{$\sum\limits^{\infty}_{m=0}$}(-1)^m{a-b \choose m} \otimes
e^m t^m=1\otimes(1-et)^{a-b}.
\end{eqnarray*}
From (\ref{fir-e2}), (\ref{fir-eeee2}) and (\ref{fom2}), we obtain
that
\\[4pt]
\hspace*{60pt}$v_au_b=\sum\limits^{\infty}_{r,s=0}{\displaystyle\frac{(-1)^s}{r!s!}}h^{[r]}_a
e^r h^{[s]}_{-b} e^s
t^{r+s}$\\[4pt]\hspace*{60pt}$\phantom{v_au_b}=\sum\limits^{\infty}_{m=0}\sum\limits_{r+s=m}{\displaystyle\frac{(-1)^s}{r!s!}}
h^{[r]}_{a}h^{[s]}_{-b-r}e^m t^m$
\\[4pt]
\hspace*{60pt}$\phantom{v_au_b}=\sum\limits^{\infty}_{m=0}{\displaystyle{a+b+m-1\choose
m}}e^m t^m=(1-et)^{-(a+b)}.$\QED
\begin{cor}\rm \label{coro}  For any $a\in\mathbb{F}$, the elements $F_a$ and $u_a$
are invertible with $F^{-1}_a=\mathcal {F}_a$, $u^{-1}_a=v_{-a}$. In
particular, $F^{-1}=\mathcal{F}$, $u^{-1}=v$.
\end{cor}
\begin{lem}\rm \label{lemm3}
For any $a\in\mathbb{F}$ and $r\in\mathbb{Z}_+$, one has $
\Delta_0(h^{[r]})=\mbox{$\sum\limits^r_{i=0}$}{r\choose
i}h^{[i]}_{-a}\otimes h^{[r-i]}_{a}. $
 In particular, one has $\Delta_0(h^{[r]})=\mbox{$\sum\limits^r_{i=0}$}{r\choose i}h^{[i]}\otimes
 h^{[r-i]}$.
\end{lem}
\noindent{\it Proof.}\ \ It can be proved by induction on $r$.\QED
\begin{lem}\rm \label{lemm4}
The element $\mathcal
{F}=\mbox{$\sum\limits^{\infty}_{r=0}$}\frac{(-1)^r}{r!}h^{[r]}\otimes
e^rt^r$ is a Drinfel'd twist on $U(\mathcal {L})[[t]]$ .
\end{lem}
\noindent{\it Proof.}\ \ It can be proved directly by using the
similar arguments as those presented in the proof of
\cite[Proposition 2.5]{HW}.\QED\vskip4pt\par Now we can perform the
process of  twisting the standard Hopf structure $(U(\mathcal
{L}),m,\iota,\Delta_0, S_0,\epsilon)$ defined in (\ref{hop}) by the
Drinfel'd twist $\mathcal {F}$ constructed above. The following
lemmas are very useful to our main results.
\begin{lem}\rm \label{lemm6}
For $p,q\in\frac{1}{2}+\mathbb{Z}$ and  $s\in\mathbb{Z}_+$, one has
\begin{eqnarray*} Y_pY_q^s=Y_q^sY_p-s(p-q)Y_q^{s-1}M_{p+q}. \end{eqnarray*}
\end{lem}
\noindent{\it Proof.}\ \ It follows from Lemma \ref{Legr4} that
\begin{eqnarray}\label{yy}
Y_pY_q^s=\mbox{$\sum\limits^s_{i=0}$}(-1)^i{s\choose
i}Y_q^{s-i}({\rm ad\,} Y_q)^{i}Y_p.
\end{eqnarray}
By (\ref{LB}), we have
\begin{eqnarray}\label{yyy}
({\rm ad\,}Y_q)^{i}Y_p=\left\{
\begin{array}{ll}
Y_p, &  i=0,\\
(p-q)M_{p+q}, &  i=1,\\
 0,&{\rm others}.
\end{array}
\right.
\end{eqnarray}
Substituting (\ref{yyy}) in (\ref{yy}) we get the result.\QED
\begin{lem}\rm \label{lemm7}   For $a\in\mathbb{F}$, $n,k\in\mathbb{Z}$, $p\in\frac{1}{2}+\mathbb{Z}$, $h=\frac{1}{n_0}L_0$ and $e=M_{n_0}$, we have
\begin{eqnarray*}
&&(L_n\otimes 1)F_a=F_{a-\frac{n}{n_0}}(L_n\otimes 1),\\[3pt]
&&(M_k\otimes 1)F_a=F_{a-\frac{k}{n_0}}(M_k\otimes1),\\[3pt]
&&(Y_p\otimes 1)F_a=F_{a-\frac{p}{n_0}}(Y_p\otimes 1),\\[3pt]
&&(1\otimes M_k)F_a=F_a(1\otimes M_k), (1\otimes
Y_p)F_a=F_a(1\otimes
Y_P),\\[3pt]&&(1\otimes L_n)F_a=F_a(1\otimes L_n)+n_0F_{a+1}(h^{(1)}_a\otimes
M_{n+n_0}t).
\end{eqnarray*}
\end{lem}
\noindent{\it Proof.}\ \ The former three equations can be directly
obtained by the definition of $F_a$ and Lemma \ref{lemm1}(i). Since
both $Y_p$ and $M_k$ commute with $e$, the next two become obvious.
It is left to verify the last one. From equations (\ref{fir-e}),
(\ref{fom1}) and Lemma \ref{Legr4}, one has
\\[4pt]\hspace*{4ex}$\dis
(1\otimes
L_n)F_a=\mbox{$\sum\limits^{\infty}_{r=0}$}\frac{1}{r!}h^{(r)}_a\otimes
L_n
e^r t^r=
\mbox{$\sum\limits^{\infty}_{r=0}$}\frac{1}{r!}h^{(r)}_a\otimes\big(
 \mbox{$\sum\limits^{r}_{i=0}$}(-1)^i{r\choose i}e^{r-i}({\rm ad\,}
e)^i
L_n\big)t^r$\\[4pt]\hspace*{4ex}$\dis\phantom{(1\otimes
L_n)F_a}
=\mbox{$\sum\limits^{\infty}_{i=0}\sum\limits^{\infty}_{r=0}$}\frac{(-1)^i}{(r+i)!}{r+i\choose
i}h^{(r+i)}_a\otimes e^r({\rm ad\,} e)^iL_n t^{r+i} $\\
[4pt]\hspace*{4ex}$\dis\phantom{(1\otimes L_n)F_a} =
\mbox{$\sum\limits^{\infty}_{i=0}\sum\limits^{\infty}_{r=0}$}\frac{(-1)^i}{r!i!}h^{(r)}_{a+i}h^{(i)}_{a}\otimes
e^r({\rm ad\,} e)^iL_n t^{r+i}$\\[4pt]\hspace*{4ex}$\dis\phantom{(1\otimes L_n)F_a}
=\mbox{$\sum\limits^{\infty}_{i=0}$}\frac{(-1)^i}{i!}F_{a+i}\big(h^{(i)}_a\otimes({\rm
ad\,} e)^i L_n t^{i}\big)$ \\[4pt]\hspace*{4ex}$\dis\phantom{(1\otimes
L_n)F_a}=F_a(1\otimes L_n)+n_0F_{a+1}\big(h^{(1)}_a\otimes
M_{n+n_0}t\big). $ \QED
\begin{lem}\rm \label{lemm8}
 For $a\in\mathbb{F}$, $n,k\in\mathbb{Z}$, $p\in\frac{1}{2}+\mathbb{Z}$, $h=\frac{2}{n_0}L_0$ and $e=Y_{\frac{n_0}{2}}$, we have
\begin{eqnarray*}
(L_n\otimes 1)F_a&=&F_{a-\frac{2n}{n_0}}(L_n\otimes 1),\ \ \
(M_k\otimes 1)F_a=F_{a-\frac{2k}{n_0}}(M_k\otimes1), \\
(Y_p\otimes 1)F_a&=&F_{a-\frac{2p}{n_0}}(Y_p\otimes 1),\ \ \
(1\otimes M_k)F_a=F_a(1\otimes M_k),\\
(1\otimes L_n)F_a &=&F_a(1\otimes
L_n)-{\displaystyle\frac{n-n_0}{2}}F_{a+1}\big(h^{(1)}_a\otimes
Y_{n+\frac{n_0}{2}}\big)t\\&
+&{\displaystyle\frac{n(n-n_0)}{4}}F_{a+2}\big(h^{(2)}_{a}\otimes
M_{n+n_0}\big)t^2,\hspace*{-100ex}\\
(1\otimes Y_p)F_a&=&F_a(1\otimes
Y_p)-(p-{\displaystyle\frac{n_0}{2}})F_{a+1}(h^{(1)}_a\otimes
M_{p+\frac{n_0}{2}})t.\hspace*{-140ex}
\end{eqnarray*}
\end{lem}
\noindent{\it Proof.}\ \ It only needs to verify the last two
formulas since the other four are obvious because of Lemma
\ref{lemm1}(ii). By (\ref{fir-e}), (\ref{fom1}) and Lemma
\ref{Legr4}, one has
\\[4pt]\hspace*{4ex}$\dis
 (1\otimes L_n)F_a=\mbox{$\sum\limits^{\infty}_{r=0}$}\frac{1}{r!}h^{(r)}_a\otimes L_n e^r
 t^r$
\\[4pt]\hspace*{4ex}$\dis\phantom{(1\otimes L_n)F_a}
 \mbox{$\sum\limits^{\infty}_{r=0}$}\frac{1}{r!}h_a^{(r)}\otimes\big(\mbox{$\sum\limits^{r}_{i=0}$}(-1)^i{r\choose
 i}e^{r-i}({\rm ad\,} e)^iL_n\big)t^r$\\[4pt]\hspace*{4ex}$\dis\phantom{(1\otimes L_n)F_a}
 =\mbox{$\sum\limits^{\infty}_{i=0}\sum\limits^{\infty}_{r=0}$}\frac{(-1)^i}
 {(r+i)!}{r+i\choose i}h^{(r+i)}_{a}\otimes e^r ({\rm ad\,} e)^iL_n t^{r+i}$\\[4pt]\hspace*{4ex}$\dis\phantom{(1\otimes L_n)F_a}
 =\mbox{$\sum\limits^{\infty}_{i=0}\sum\limits^{\infty}_{r=0}$}\frac{(-1)^i}{r!i!}h^{(r)}_{a+i}h^{(i)}_a\otimes e^r({\rm ad\,} e)^i
 L_nt^{r+i}$\\[4pt]\hspace*{4ex}$\dis\phantom{(1\otimes L_n)F_a}
 =\mbox{$\sum\limits^{\infty}_{i=0}$}\frac{(-1)^i}{i!}F_{a+i}\big(h^{(i)}_a\otimes({\rm
 ad\,} e)^iL_n\big)t^i$\\[4pt]\hspace*{4ex}$\dis\phantom{(1\otimes L_n)F_a}
 =F_a(1\otimes L_n)-\frac{n-n_0}{2}F_{a+1}\big(h^{(1)}_a\otimes
 Y_{n+\frac{n_0}{2}}\big)t$\\[4pt]\hspace*{4ex}$\dis\phantom{(1\otimes L_n)F_a}+\frac{n(n-n_0)}{4}F_{a+2}\big(h^{(2)}_{a}\otimes
 M_{n+n_0}\big)t^2.$\\
Again by (\ref{fir-e}), (\ref{fom1}) and Lemma \ref{lemm6}, one has
\\[4pt]\hspace*{4ex}$\dis
(1\otimes Y_p)F_a=
\mbox{$\sum\limits^{\infty}_{r=0}$}\frac{1}{r!}h^{(r)}_{a}\otimes\big(e^rY_p-r(p-\frac{n_0}{2})e^{r-1}
M_{p+\frac{n_0}{2}}\big)t^r$\\[4pt]\hspace*{4ex}$\dis\phantom{(1\otimes Y_p)F_a}=
F_a(1\otimes
Y_p)-(p-\frac{n_0}{2})\mbox{$\sum\limits^{\infty}_{r=1}$}\frac{1}{(r-1)!}h^{(r)}_a\otimes
e^{r-1} M_{p+\frac{n_0}{2}}t^r$\\[4pt]\hspace*{4ex}$\dis\phantom{(1\otimes Y_p)F_a}=
F_a(1\otimes
Y_p)-(p-\frac{n_0}{2})\mbox{$\sum\limits^{\infty}_{r=0}$}\frac{1}{r!}h^{(r+1)}_a\otimes
e^{r}M_{p+\frac{n_0}{2}}t^{r+1}$\\[4pt]\hspace*{4ex}$\dis\phantom{(1\otimes Y_p)F_a}=
F_a(1\otimes
Y_p)-(p-\frac{n_0}{2})\mbox{$\sum\limits^{\infty}_{r=0}$}\frac{1}{r!}h^{(r)}_{a+1}h^{(1)}_a\otimes
e^{r}M_{p+\frac{n_0}{2}}t^{r+1}$\\[4pt]\hspace*{4ex}$\dis\phantom{(1\otimes Y_p)F_a}=
F_a(1\otimes Y_p)-(p-\frac{n_0}{2})F_{a+1}(h^{(1)}_a\otimes
M_{p+\frac{n_0}{2}})t. $
\QED
\begin{lem}\rm \label{lemm9}
 For $a\in\mathbb{F}$, $n,k\in\mathbb{Z}$, $p\in\frac{1}{2}+\mathbb{Z}$, $h=\frac{1}{n_0}L_0$ and $e=M_{n_0}$, we have
 \begin{eqnarray}
 &&L_n
 u_a=u_{a+\frac{n}{n_0}}\big(L_n-n_0M_{n+n_0}h^{[1]}_{1-a}t\big),\label{fom3}\\
&&
 M_ku_a=u_{a+\frac{k}{n_0}}M_k,\ \
 Y_pu_a=u_{a+\frac{p}{n_0}}Y_p. \label{fom333}
 \end{eqnarray}
\end{lem}
\noindent{\it Proof.}\ \ From Lemma \ref{lemm1}(i), and since both
$Y_p$ and $M_k$ commute with
 $e$, one can easily obtain the last two equations of (\ref{fom333}). By (\ref{fir-e}), (\ref{fom2}) and Lemma \ref{Legr4},
 we have
\\[4pt]\hspace*{4ex}$\dis
 L_nu_a=\mbox{$\sum\limits^{\infty}_{r=0}$}\frac{(-1)^r}{r!}h^{[r]}_{-a-\frac{n}{n_0}}L_ne^rt^r$\\
[4pt]\hspace*{4ex}$\dis\phantom{L_nu_a} =\mbox{$\sum\limits^{\infty}_{r=0}$}\frac{(-1)^r}{r!}h^{[r]}_{-a-\frac{n}{n_0}}(e^rL_n+rn_0e^{r-1}M_{n+n_0})t^r$\\
[4pt]\hspace*{4ex}$\dis\phantom{L_nu_a}=u_{a+\frac{n}{n_0}}L_n+n_0\mbox{$\sum\limits^{\infty}_{r=1}$}\frac{(-1)^r}{(r-1)!}h^{[r]}_{-a-\frac{n}{n_0}}e^{r-1}M_{n+n_0}t^{r}
$\\[4pt]\hspace*{4ex}$\dis\phantom{L_nu_a} =u_{a+\frac{n}{n_0}}L_n-n_0\mbox{$\sum\limits^{\infty}_{r=0}$}\frac{(-1)^r}{r!}h^{[r+1]}_{-a-\frac{n}{n_0}}e^{r}M_{n+n_0}t^{r+1}$\\
[4pt]\hspace*{4ex}$\dis\phantom{L_nu_a}
=u_{a+\frac{n}{n_0}}L_n-n_0\mbox{$\sum\limits^{\infty}_{r=0}$}\frac{(-1)^r}{r!}h^{[r]}_{-a-\frac{n}{n_0}}h^{[1]}_{-a-\frac{n}{n_0}-r}e^{r}M_{n+n_0}t^{r+1}$\\
[4pt]\hspace*{4ex}$\dis\phantom{L_nu_a}
=u_{a+\frac{n}{n_0}}L_n-n_0\mbox{$\sum\limits^{\infty}_{r=0}$}\frac{(-1)^r}{r!}h^{[r]}_{-a-\frac{n}{n_0}}e^rh^{[1]}_{-a-\frac{n}{n_0}}M_{n+n_0}t^{r+1}$\\
[4pt]\hspace*{4ex}$\dis\phantom{L_nu_a}
 =u_{a+\frac{n}{n_0}}\big(L_n-n_0h^{[1]}_{-a-\frac{n}{n_0}}M_{n+n_0}t\big)
 =u_{a+\frac{n}{n_0}}\big(L_n-n_0M_{n+n_0}h^{[1]}_{1-a}t\big).$
 \QED
\begin{lem}\rm \label{lemm10} For $a\in\mathbb{F}$, $n,k\in\mathbb{Z}$, $p\in\frac{1}{2}+\mathbb{Z}$, $h=\frac{2}{n_0}L_0$ and $e=Y_{\frac{n_0}{2}}$, we have
\begin{eqnarray}
&&M_ku_a=u_{a+\frac{2k}{n_0}}M_k,\ \ \ \ \ \
Y_pu_a=u_{a+\frac{2p}{n_0}}\big(Y_p+(p-\frac{n_0}{2})M_{p+\frac{n_0}{2}}h^{[1]}_{1-a}t\big),\ \ \ \ \ \ \ \ \label{fom8}\\[4pt]
&&L_nu_a=u_{a+\frac{2n}{n_0}}\big(L_n+\frac{n-n_0}{2}Y_{n+\frac{n_0}{2}}h^{[1]}_{1-a}t+\frac{n(n-n_0)}{4}M_{n+n_0}h^{[2]}_{2-a}t^2\big).\
\ \ \ \ \ \ \ \label{fom9}
\end{eqnarray}
\end{lem}
 \noindent{\it Proof.}\ \  The first equaton of (\ref{fom8}) is obvious by Lemma  \ref{lemm1}(ii) and since $M_k$ commutes with
 $e$. From (\ref{fir-e}), (\ref{fom2}), Lemmas \ref{Legr4} and \ref{lemm1}(ii), one has
\\[4pt]\hspace*{4ex}$\dis
 L_nu_a=\mbox{$\sum\limits^{\infty}_{r=0}$}\frac{(-1)^r}{r!}h^{[r]}_{-a-\frac{2n}{n_0}}L_ne^rt^r
$\\[4pt]\hspace*{4ex}$\dis\phantom{L_nu_a}=\mbox{$\sum\limits^{\infty}_{r=0}$}\frac{(-1)^r}{r!}h^{[r]}_{-a-\frac{2n}{n_0}}\big(
\mbox{$\sum\limits^{r}_{i=0}$}(-1)^i{r\choose
 i}e^{r-i}({\rm ad\,} e)^iL_n\big)t^r$\\
[4pt]\hspace*{4ex}$\dis\phantom{L_nu_a}
=\mbox{$\sum\limits^{\infty}_{i=0}$}\mbox{$\sum\limits^{\infty}_{r=0}$}\frac{(-1)^r}{(r+i)!}{r+i\choose
i}h^{[r+i]}_{-a-\frac{2n}{n_0}}e^r({\rm ad\,} e)^iL_nt^{r+i}$\\
[4pt]\hspace*{4ex}$\dis\phantom{L_nu_a}=\mbox{$\sum\limits^{\infty}_{i=0}$}\mbox{$\sum\limits^{\infty}_{r=0}$}\frac{(-1)^r}{r!i!}h^{[r]}_{-a-\frac{2n}{n_0}}h^{[i]}_{-a-\frac{2n}{n_0}-r}e^r({\rm
ad\,}e)^i L_nt^{r+i} $\\[4pt]
\hspace*{4ex}$\dis\phantom{L_nu_a}
=\mbox{$\sum\limits^{\infty}_{i=0}$}\mbox{$\sum\limits^{\infty}_{r=0}$}\frac
{(-1)^r}{r!i!}h^{[r]}_{-a-\frac{2n}{n_0}}e^rh^{[i]}_{-a-\frac{2n}{n_0}}({\rm
ad\,} e)^iL_nt^{r+i} $\\[4pt]\hspace*{4ex}$\dis\phantom{L_nu_a}=u_{a+\frac{2n}{n_0}}\big(
\mbox{$\sum\limits^{\infty}_{i=0}$}\frac{1}{i!}h^{[i]}_{-a-\frac{2n}{n_0}}({\rm
ad\,}e)^iL_nt^i\big)$\\
[4pt]\hspace*{4ex}$\dis\phantom{L_nu_a}
=u_{a+\frac{2n}{n_0}}\big(L_n+\frac{n-n_0}{2}h^{[1]}_{-a-\frac{2n}{n_0}}Y_{n+\frac{n_0}{2}}t+\frac{n(n-n_0)}{4}h^{[2]}_{-a-\frac{2n}{n_0}}M_{n+n_0}t^2\big)$\\
[4pt]\hspace*{4ex}$\dis\phantom{L_nu_a}
=u_{a+\frac{2n}{n_0}}\big(L_n+\frac{n-n_0}{2}Y_{n+\frac{n_0}{2}}h^{[1]}_{1-a}t+\frac{n(n-n_0)}{4}M_{n+n_0}h^{[2]}_{2-a}t^2\big).
$\\ In addition, by Lemma \ref{lemm6}, one has
\\[4pt]\hspace*{4ex}$\dis
Y_pu_a=\mbox{$\sum\limits^{\infty}_{r=0}$}\frac{(-1)^r}{r!}h^{[r]}_{-a-\frac{2p}{n_0}}Y_pe^rt^r
$\\[4pt]\hspace*{4ex}$\dis\phantom{Y_pu_a}=\mbox{$\sum\limits^{\infty}_{r=0}$}\frac{(-1)^r}{r!}h^{[r]}_{-a-\frac{2p}{n_0}}\big(e^rY_p-r(p-\frac{n_0}{2})e^{r-1}M_{p+\frac{n_0}{2}}\big)t^r$\\
[4pt]\hspace*{4ex}$\dis\phantom{Y_pu_a}
=u_{a+\frac{2p}{n_0}}Y_p-(p-\frac{n_0}{2})\mbox{$\sum\limits^{\infty}_{r=1}$}\frac{(-1)^r}{(r-1)!}h^{[r]}_{-a-\frac{2p}{n_0}}e^{r-1}M_{p+\frac{n_0}{2}}t^r$\\
[4pt]\hspace*{4ex}$\dis\phantom{Y_pu_a}
=u_{a+\frac{2p}{n_0}}Y_p+(p-\frac{n_0}{2})\mbox{$\sum\limits^{\infty}_{r=0}$}\frac{(-1)^r}{r!}h^{[r]}_{-a-\frac{2p}{n_0}}h^{[1]}_{-a-\frac{2p}{n_0}-r}e^rM_{p+\frac{n_0}{2}}t^{r+1}$\\
[4pt]\hspace*{4ex}$\dis\phantom{Y_pu_a}
=u_{a+\frac{2p}{n_0}}Y_p+(p-\frac{n_0}{2})\mbox{$\sum\limits^{\infty}_{r=0}$}\frac{(-1)^r}{r!}h^{[r]}_{-a-\frac{2p}{n_0}}e^rh^{[1]}_{-a-\frac{2p}{n_0}}M_{p+\frac{n_0}{2}}t^{r+1}$\\
[4pt]\hspace*{4ex}$\dis\phantom{Y_pu_a}
=u_{a+\frac{2p}{n_0}}Y_p+(p-\frac{n_0}{2})u_{a+\frac{2p}{n_0}}(h^{[1]}_{-a-\frac{2p}{n_0}}M_{p+\frac{n_0}{2}}t)$\\[4pt]\hspace*{4ex}$\dis\phantom{Y_pu_a}
=u_{a+\frac{2p}{n_0}}\big(Y_P+(p-\frac{n_0}{2})M_{p+\frac{n_0}{2}}h^{[1]}_{1-a}t\big).$
 \QED\vspace{4mm}\par
Now we have enough in hand to prove our main results in this
paper.\\[4mm]
 \noindent{\it Proof of Theorem \ref{theo1}}~~\rm\  By Lemmas
\ref{Legr}, \ref{lemm2}, Corollary \ref{coro}, Lemmas \ref{lemm4} and \ref{lemm7}, we have\\
[4pt]\hspace*{4ex}$\dis
 \Delta(L_n)=\mathcal
{F}\cdot \Delta_0(L_n)\cdot \mathcal{F}^{-1}=\mathcal {F}\cdot
(L_n\otimes
1)\cdot F+\mathcal {F}\cdot (1\otimes L_n)\cdot F$\\
[4pt]\hspace*{4ex}$\dis\phantom{\Delta(L_n)} =\mathcal {F}\cdot
F_{-\frac{n}{n_0}}\cdot L_n\otimes 1+\mathcal{F}\cdot(F\cdot
1\otimes L_n+n_0F_1\cdot h^{(1)}\otimes
M_{n_0+n}t)$\\
[4pt]\hspace*{4ex}$\dis\phantom{\Delta(L_n)}
=1\otimes(1-et)^{\frac{n}{n_0}}\cdot L_n\otimes 1+1\otimes
L_n+n_0\otimes(1-et)^{-1}\cdot (h^{(1)}\otimes
M_{n+n_0}t)$\\
[4pt]\hspace*{4ex}$\dis\phantom{\Delta(L_n)} =1\otimes L_n+
L_n\otimes (1-et)^{\frac{n}{n_0}}+n_0h\otimes (1-et)^{-1}M_{n+n_0}t.
$\\
[4pt]\hspace*{4ex}$\dis\Delta(M_k)=\mathcal{F}\cdot (M_k\otimes
1)\cdot F+\mathcal{F}\cdot (1\otimes M_k)\cdot F$\\
[4pt]\hspace*{4ex}$\dis\phantom{\Delta(M_k)} =\mathcal{F}\cdot
F_{-\frac{k}{n_0}}\cdot (M_k\otimes 1)+\mathcal{F}\cdot F(1\otimes
M_k)$
\\ \hspace*{4ex}$\dis\phantom{\Delta(M_k)}
=1\otimes(1-et)^{\frac{k}{n_0}}\cdot(M_k\otimes 1)+1\otimes M_k
=M_k\otimes(1-et)^{\frac{k}{n_0}}+1\otimes M_k. $\\
[4pt]\hspace*{4ex}$\dis \Delta(Y_p)=\mathcal
{F}\Delta_0(Y_p)\mathcal {F}^{-1}=\mathcal{F}\cdot Y_p\otimes 1\cdot
F+\mathcal{F}\cdot 1\otimes Y_p\cdot F$\\
[4pt]\hspace*{4ex}$\dis\phantom{\Delta(Y_p)} =\mathcal{F}\cdot
F_{-\frac{p}{n_0}}\cdot Y_p\otimes 1+\mathcal{F}\cdot F\cdot
1\otimes
Y_p$\\
\hspace*{4ex}$\dis\phantom{\Delta(Y_p)} =1\otimes
(1-et)^{\frac{p}{n_0}}\cdot Y_p\otimes 1+1\otimes Y_p
=Y_p\otimes(1-et)^{\frac{p}{n_0}}+1\otimes Y_p. $\\
Again by Lemmas \ref{Legr}, \ref{lemm2}, Corollary \ref{coro} and Lemma \ref{lemm9}, we have\\[4pt]\hspace*{4ex}$\dis
S(L_n)=-vL_nu =-v\cdot
u_{\frac{n}{n_0}}(L_n-n_0M_{n+n_0}h^{[1]}_{1}t)$\\[4pt]\hspace*{4ex}$\dis\phantom{
S(L_n)}
=-(1-et)^{-\frac{n}{n_0}}(L_n-n_0M_{n+n_0}h^{[1]}_{1}t).$\\[4pt]\hspace*{4ex}$\dis S(M_k)=
u^{-1}S_0(M_k)u=-v\cdot M_k\cdot u =-vu_{\frac{k}{n_0}}M_k
=-(1-et)^{-\frac{k}{n_0}}\cdot M_k.$\\[4pt]\hspace*{4ex}$\dis
S(Y_p)=u^{-1}S_0(Y_p)u=-v\cdot Y_p\cdot u=-v\cdot
u_{\frac{p}{n_0}}\cdot Y_p =-(1-et)^{-\frac{p}{n_0}}\cdot Y_p. $\\
Hence, we get the results. \QED

\vspace{4mm}

\noindent{\it Proof of Theorem \ref{theo2}}~~\rm  By Lemmas
\ref{Legr}, \ref{lemm2}, Corollary \ref{coro}, Lemmas \ref{lemm4}
and \ref{lemm8}, we have\\
[4pt]\hspace*{4ex}$\dis \Delta(L_n)=\mathcal
{F}\Delta_0(L_n)\mathcal {F}^{-1}=\mathcal {F}\cdot L_n\otimes
1\cdot F+\mathcal {F}\cdot 1\otimes L_n\cdot
F$\\[4pt]\hspace*{4ex}$\dis\phantom{\Delta(L_n)}
=\mathcal {F}\cdot F_{-\frac{2n}{n_0}}\cdot L_n\otimes 1+\mathcal
{F}\cdot \big(F\cdot 1\otimes L_n+\frac{n_0-n}{2}F_1\cdot h\otimes
Y_{n+\frac{n_0}{2}}t$\\[4pt]\hspace*{4ex}$\dis\phantom{\Delta(L_n)}\ \ +\frac{n(n-n_0)}{4}F_2\cdot h^{(2)}\otimes
M_{n+n_0}t^2\big)$\\[4pt]\hspace*{4ex}$\dis\phantom{\Delta(L_n)}
=1\otimes(1-et)^{\frac{2n}{n_0}}L_n\otimes 1+1\otimes
L_n+\frac{n_0-n}{2}\otimes (1-et)^{-1}h\otimes
Y_{n+\frac{n_0}{2}}t$\\[4pt]\hspace*{4ex}$\dis\phantom{\Delta(L_n)}\ \ +\frac{n(n-n_0)}{4}\otimes
(1-et)^{-2}\cdot h^{(2)}\otimes
M_{n+n_0}t^2$\\[4pt]\hspace*{4ex}$\dis\phantom{\Delta(L_n)}
=L_n\otimes (1-et)^{\frac{2n}{n_0}}+1\otimes L_n
+\frac{n_0-n}{2}h\otimes (1-et)^{-1}Y_{n+\frac{n_0}{2}}t$
\\[4pt]\hspace*{4ex}$\dis\phantom{\Delta(L_n)}\ \ \ +\frac{n(n-n_0)}{4}h^{(2)}\otimes
(1-et)^{-2}M_{n+n_0}t^2$.\\

\begin{eqnarray*}
 \Delta(M_k)&=&\mathcal {F}\cdot
M_k\otimes 1\cdot F+\mathcal {F}\cdot 1\otimes M_k \cdot F\\
&=&\mathcal {F}\cdot F_{-\frac{2k}{n_0}}\cdot M_k\otimes 1+\mathcal
{F}\cdot F\cdot 1\otimes M_k\\[-6pt]
&=&1\otimes (1-et)^{\frac{2k}{n_0}}\cdot M_k\otimes 1+1\otimes M_k
=M_k\otimes (1-et)^{\frac{2k}{n_0}}+1\otimes M_k.\\[4pt]
\Delta(Y_p)&=&\mathcal {F}\Delta_0(L_n)\mathcal {F}^{-1}=\mathcal {F}\cdot Y_p\otimes 1\cdot F+\mathcal{F}\cdot 1\otimes Y_p \cdot F\\
&=&\mathcal {F} F_{-\frac{2p}{n_0}}\cdot Y_P\otimes 1+\mathcal
{F}\cdot\big(F(1\otimes Y_p)-(p-\frac{n_0}{2})F_1(h\otimes
M_{p+\frac{n_0}{2}})t\big)\\
&=&1\otimes (1-et)^{\frac{2p}{n_0}}\cdot Y_p\otimes 1+1\otimes
Y_p-(p-\frac{n_0}{2})\otimes (1-et)^{-1}\cdot h\otimes
M_{p+\frac{n_0}{2}}t\\
&=&1\otimes Y_p+Y_p\otimes
(1-et)^{\frac{2p}{n_0}}-(p-\frac{n_0}{2})h\otimes
(1-et)^{-1}M_{p+\frac{n_0}{2}} t.
\end{eqnarray*}
Again by Lemmas \ref{Legr}, \ref{lemm2}, Corollary \ref{coro} and
Lemma \ref{lemm10}, we have
\begin{eqnarray*}
S(L_n)&=&u^{-1}S_0(L_n)u=-vL_nu\\[4pt]
&=&-vu_{\frac{2n}{n_0}}(L_n+\frac{n-n_0}{2}Y_{n+\frac{n_0}{2}}h^{[1]}_1
t+\frac{n(n-n_0)}{4}M_{n+n_0}h^{[2]}_{2}t^2)\\[4pt]
&=&-(1-et)^{-\frac{2n}{n_0}}\big(L_n+\frac{n-n_0}{2}Y_{n+\frac{n_0}{2}}h^{[1]}_1
t+\frac{n(n-n_0)}{4}M_{n+n_0}h^{[2]}_{2}t^2\big).\\[4pt]
S(Y_p)&=&u^{-1}S_0(Y_p)u=-vY_pu =-v\cdot
u_{\frac{2p}{n_0}}\cdot\big(Y_p+(p-\frac{n_0}{2})M_{p+\frac{n_0}{2}}h^{[1]}_1t\big)\\[4pt]
&=&-(1-et)^{-\frac{2p}{n_0}}\big(Y_p+(p-\frac{n_0}{2})M_{p+\frac{n_0}{2}}h^{[1]}_1t\big).\\[4pt]
S(M_k)&=&u^{-1}S_0(M_k)u=-v\cdot M_ku=-v\cdot
u_{\frac{2k}{n_0}}M_k=-(1-et)^{-\frac{2k}{n_0}}M_k.
\end{eqnarray*}
So the proof is complete!\QED
\vspace{2mm}\\
\begin{rem}\rm In this paper, we have presented two
kinds of Hopf algebraic structures on $U(\LL)[[t]]$ using the
Drinfel'd twists. During the process of constructing the Drinfel'd
twists, we see that any of them is definitely determined by some
classical Yang-Baxter $r$-matrix $r$ (namely, the Lie bialgebra
structures of $\LL$). So any two different elements $h,e\in\LL$ such
that $[h,e]=e$ can determine a Drinfel'd twist and thus a Hopf
algebraic structure. This is one of the reasons why it is difficult
to determine all Hopf algebraic structures on $U(\LL)[[t]]$. It is
sure that there exist other Hopf algebraic structures different from
that given in our paper. One clear example is to take
$h=\frac{L_0}{n_0}$ and $e=L_{n_0}$ for a fixed nonzero integer
$n_0$. It is easy to see $[h,e]=e$. Thus one can get another Hopf
algebraic structure using the similar arguments as above.
\end{rem}

\vspace{4mm} \noindent\bf{\footnotesize Acknowledgements}\quad\rm
 {\footnotesize
This research was supported by National Natural Science Foundation
grants of
China (10825101).}\\[4mm]

\noindent{\bbb{References}}

\end{document}